# "A PRIORI" ESTIMATES FOR SOLUTIONS OF SINGULAR ELLIPTIC INEQUALITIES AND APPLICATIONS

STEFANO PIGOLA, MARCO RIGOLI, AND ALBERTO G. SETTI

*Dedicated to the memory of Franca Burrone Rigoli*

ABSTRACT. We obtain a maximum principle, and "a priori" upper estimates for solutions of a class of non linear singular elliptic differential inequalities on Riemannian manifolds under the sole geometrical assumption of volume growth conditions. Various applications of the results obtained are presented.

## 0. Introduction

Let $(M, <, >)$ be a smooth, connected, non-compact, complete Riemannian manifold of dimension $m \geq 2$. We fix an origin $o$, and denote by $r(x)$ the distance function from $o$, and by $B_t = \{x \in M : r(x) < t\}$ and $\partial B_t = \{x \in M : r(x) = t\}$ the geodesic ball and sphere of radius $t > 0$ centered at $o$, respectively.

Let $\varphi$ be a real valued function in $C^1\big((0, +\infty)\big) \cap C^o\big([0, +\infty)\big)$ satisfying the following structural conditions:

(0.1)   i) $\varphi(0) = 0$;   ii) $\varphi(t) > 0 \ \forall t > 0$;   iii) $\varphi(t) \leq At^\delta \ \forall t \geq 0$,

for some positive constants $A$ and $\delta$. We will consider the differential operator defined, for $u \in C^1(M)$, by

(0.2) $$\mathrm{div}\big(|\nabla u|^{-1} \varphi(|\nabla u|) \nabla u\big),$$

and which we refer to as the $\varphi$-Laplacian. It is understood that, if the the vector field in brackets is not $C^1$, then the divergence in (0.2) must be considered in distributional sense, and equations and inequalities will be interpreted accordingly.

We remark that the $\varphi$-Laplacian arises naturally when considering the Euler–Lagrange equation associated to the energy functional

$$\Lambda(u) = \int \Phi(|\nabla u|),$$

where $\Phi(t) = \int_0^t \varphi(s)\,ds$.

As important natural examples we mention

1. the Laplace–Beltrami operator, $\Delta u$, corresponding to $\varphi(t) = t$;
2. or, more generally, the $p$-Laplacian, $\mathrm{div}(|\nabla u|^{p-2}\nabla u)$, $p > 1$, corresponding to $\varphi(t) = t^{p-1}$;
3. the generalized mean curvature operator, $\mathrm{div}\big(\frac{\nabla u}{(1+|\nabla u|^2)^\alpha}\big)$, $\alpha > 0$, corresponding to $\varphi(t) = t/(1+t^2)^\alpha$.

In the Euclidean setting the $\varphi$-Laplacian has been studied by many authors. We mention in particular the seminal papers of J. Serrin [S1], and Serrin and L.A. Peletier, [PeS], where Liouville type theorems are established. See also [S2], [PuSZ], [PuS], where the maximum principle for the $\varphi$-Laplacian is investigated, and [FLS] devoted to the existence and uniqueness of radial ground states.

It is often vital, in solving a variety of problems, to be able to establish suitable "a-priori" estimates for solutions of equations or differential inequalities of the form

(0.3) $$\mathrm{div}\big(|\nabla u|^{-1}\varphi(|\nabla u|)\nabla u\big) \geq F(x,u).$$

In geometrical problems, such as those presented in Section 3 below, this type of information frequently leads quickly to geometrical conclusions. In a Riemannian setting, the natural assumptions are those more or less directly related to curvature. For instance, the experienced reader will certainly recognize that in the study of analytic properties of the manifold the role played by the Ricci curvature can hardly be overestimated. In fact, curvature assumptions imply estimates on the volume growth of geodesic balls, while the reverse implication is far from being true.

In the present note we concentrate our attention on volume growth assumptions and show that these suffice to obtain interesting properties of solutions to (0.3).

We next describe our main results, referring to the forthcoming paper [PRS3] for the proofs and further details.

**Theorem A.** *Let $(M, <,>)$ and $\varphi$ be as above and suppose that $b(x)$ is a continuous function on $M$ satisfying*

(0.4) $$b(x) \geq \frac{1}{Q(r(x))} \quad on\ M$$

where $Q(t)$ is positive, non-decreasing and satisfies $Q(t) = o(t^{1+\delta})$ as $t \to +\infty$. Let $f \in C^0(\mathbb{R})$ and assume that $u \in C^1(M)$ satisfy $u^* = \sup_M u < +\infty$ and

$$(0.5) \qquad \mathrm{div}\bigl(|\nabla u|^{-1}\varphi(|\nabla u|)\nabla u\bigr) \geq b(x)f(u)$$

on the set

$$(0.6) \qquad \Omega_\gamma = \{x \in M : u(x) > \gamma\},$$

for some $\gamma < u^*$. If

$$(0.7) \qquad \liminf_{r \to +\infty} \frac{Q(r)\log \mathrm{vol}\, B_r}{r^{1+\delta}} < +\infty,$$

then $f(u^*) \leq 0$.

We note that assumption (0.7) may be replaced by

$$(0.8) \qquad \liminf_{r \to +\infty} \frac{Q(r)}{r^{1+\delta}} \int_{B_r} |u|^p < +\infty \quad \text{for some } p > 0.$$

We also remark that if the vector field $|\nabla u|^{-1}\varphi(|\nabla u|)\nabla u$ is of class at least $C^1$, then the conclusion of Theorem A is equivalent to the following version of the weak maximum principle, namely, that there exists a sequence $\{x_k\}$ such that

$$u(x_k) > u^* - 1/k, \qquad Q(r(x_k))\,\mathrm{div}\bigl(|\nabla u|^{-1}\varphi(|\nabla u|)\nabla u\bigr)(x_k) < 1/k.$$

When $Q(r)$ is constant, this is precisely the weak maximum principle for the $\varphi$-Laplacian established in [RS, Theorem 5.1]. As shown in [PRS2], in the case of the Laplacian this turns out to be equivalent to the stochastic completeness of the manifold $M$. In this sense, the conclusion of Theorem A may be viewed as a stronger version of stochastic completeness. Intermediate versions of stochastic completeness have been considered by A. Grigor'yan [G1], and in a different direction by Grigor'yan and M. Kelbert [GK].

For more information on stochastic completeness, parabolicity and analytic and geometric consequences thereof, we refer the reader to the comprehensive survey paper [G2].

We also note that a result closely related to Theorem A was proved by the authors under assumptions on the Ricci curvature of $M$, see Theorems A and 1.2 in [PRS1]. For comparison, we remark that the curvature condition assumed there, which is in general far more stringent that the volume growth condition considered here, allowed us to assume that $u$ be a solution of (0.5) on the smaller set

$$\widehat{\Omega}_\gamma = \{x \in M : u(x) > \gamma,\ |\nabla u|(x) < u^* - \gamma\}.$$

Our second main result is

**Theorem B.** *Let $(M, <,>)$ and $\varphi$ be as above. Let $a(x), b(x) \in C^0(M)$ where $a(x) = a_+(x) - a_-(x)$, with that $a_\pm \geq 0$, and assume that $\|a_-\|_\infty < +\infty$ and that $b(x)$ satisfy (0.4). Assume further that, for some $H > 0$,*

$$(0.9) \qquad \frac{a_-(x)}{b(x)} \leq H \qquad on \ M.$$

*Let $u \in C^1(M)$ be a non-negative solution of*

$$(0.10) \qquad \mathrm{div}\big(|\nabla u|^{-1}\varphi(|\nabla u|)\nabla u\big) \geq b(x)u^\sigma + a(x)u \qquad on \ M,$$

*with $\sigma > \max\{1, \delta\}$. If (0.7) holds, then*

$$(0.11) \qquad u(x) \leq H^{1/(\sigma-1)} \qquad on \ M.$$

**Remarks 0.1.** a) If $a_-(x) \equiv 0$, then Theorem B implies that any non-negative solution of (0.10) is identically zero.
b) Assume that $u$ is a solution on $M$ of the differential inequality obtained from (0.10) by replacing the non-linear term $b(x)u^\sigma$ with $b(x)f(u)$. If

$$\liminf_{t \to +\infty} \frac{f(t)}{t^\sigma} > 0,$$

for some $\sigma > 1$, and the remaining assumptions of Theorem B hold, then we can conclude that $u^* < +\infty$.
c) We note that, if $\varphi(t) = t$, $\sigma \geq (m+2)/(m-2)$ where $m = \dim M \geq 3$, and $u$ is positive, then we can relax the assumption that $b(x)$ be strictly positive on $M$ with with $b(x) \geq 0$ and conclude that $u* < +\infty$ (see [BRS, Lemma 3.6]). Further, if $a_- \equiv 0$, we may conclude that inequality (0.10) has no positive solution.
d) The same observation made after Theorem A shows that the volume growth assumption (0.7) may be replaced by (0.8)

We note that sharp detailed "a-priori" estimates for solutions of inequalities of the form (0.10) involving the Laplacian are given in [RRV, Proposition 4.1]. There in fact, it was possible to deal with the case where $\frac{a_-(x)}{b(x)} \leq H(r(x))$, but the main geometric assumption is expressed in terms of bounds on the Ricci curvature of the manifold. On the other hand, to the best of the authors' knowledge, the result described in the present paper are the first of their kind obtained under volume growth assumptions only.

The techniques used in the proof of Theorem A also yield the following Phragmen- Lindeloff type result.

**Theorem C.** *Let $\varphi$ and $b$ satisfy the conditions listed in the statement of Theorem A, and let $f$ be a continuous function on $\mathbb{R}$ such that $f > 0$* on $[\Gamma, +\infty)$, for some $\Gamma \in \mathbb{R}$. Let also $\Omega$ be an unbounded domain with non-empty boundary $\partial \Omega$, and assume that $u \in C^2(\Omega) \cap C^0(\overline{\Omega})$ be bounded and satisfy

$$\begin{cases} div(|\nabla u|^{-1}\varphi(|\nabla u|)\nabla u) \geq b(x)f(u) & on \ \Omega \\ u \leq \Gamma & on \ \partial \Omega. \end{cases}$$

If the volume growth condition (0.7) holds, then $u \leq \Gamma$ on $\Omega$.

## 1. Applications

In this section we will describe some applications of the results presented in Section 1.

We begin with a version of Schwarz lemma for holomorphic maps between complex manifolds of the same dimension.

**Theorem 1.1.** *Let $(M, <,>_M)$ be a complete Kahler manifold of complex dimension $m$ and scalar curvature $s(x)$. Set $r(x) = dist_M(x, o)$ for some fixed reference $o \in M$. Let $f : M \to N$ be a holomorphic map into a complex Hermitian manifold $(N, <,>_N)$ of the same dimension and with Ricci curvature bounded from above by $R_N(z)$, $\forall z \in N$. Suppose that*

$$R_N(f(x)) < 0 \ on \ M \quad and \ R_N(f(x)) \leq -\frac{1}{Q(r(x))} \ on \ M \setminus B_{R_0},$$

*for some constant $R_0 > 0$, and some non-decreasing function $Q(t)$ defined for $t \geq 0$ such that, as $t \to +\infty$, $Q(t) = o(t^2)$. Furthermore, assume that*

$$\|s_-\|_\infty < +\infty, \qquad \frac{s_-(x)}{R_N(f(x))} \geq -mH, \quad on \ M$$

*for some $H \geq 0$. Then, indicating with $dV_M$ and $dV_N$ the volume elements respectively of $M$ and $N$, we have*

$$\left\{ \frac{f^* dV_N}{dV_M} \right\}^{\frac{1}{m}} \leq H,$$

*provided*

$$\liminf_{r \to +\infty} \frac{Q(r) \log vol B_r}{r^2} < +\infty.$$

*In particular, if $H < 1$, $f$ is volume decreasing.*

The following can be consider as a further generalization of Schwarz lemma, noting that in the (real) 2-dimensional case the holomorphic maps are smooth conformal maps. Let $f : (M, <,>_M) \to (N, <,>_N)$ be a conformal immersion. Denote by $s(x)$ the scalar curvature of $(M, <,>_M)$ and by $k(x)$ that of the pull-back metric $f^* <,>_N$. To simplify the exposition we assume that $m = \dim M \geq 3$. Then, the conformality factor of $f$ is determined by setting $f^* <,>_N = u^{\frac{4}{m-2}} <,>_M$, with $u(x) > 0$ on $M$.

**Theorem 1.2.** *Let $(M, <,>)$ be a complete manifold with scalar curvature $s(x)$ satisfying*

$$i) \; -c \leq s(x) < 0 \quad ii) \; s(x) \leq -\frac{1}{Q(t)}, \; \text{for } r(x) \gg 1.$$

*Assume that*

$$(1.1) \qquad \liminf_{r \to +\infty} \frac{Q(r) \log \text{vol} B_r}{r^2} < +\infty$$

*where $Q(t)$ is a non decreasing function defined for $t \geq 0$ satisfying $Q(t) = o(t^2)$ as $t \to +\infty$. (1.1) holds. Then, any conformal diffeomorphism of $(M, <,>)$ into itself which preserves the scalar curvature is an isometry.*

Applying Remark 0.1 c) one obtains the following nonexistence result for positive solutions of the non-compact Yamabe equation, under volume growth conditions.

**Theorem 1.3.** *Let $(M, <,>)$ be a complete Riemannian manifold with scalar curvature $s(x) \geq 0$, and let $K(x)$ be a smooth function on $M$ satisfying $K(x) \leq 0$ and $K(x) \leq -1/Q(r(x))$ for $r(x) \gg 1$, where $Q(t)$ is an increasing function such that $Q(t) = o(t^2)$ as $t \to +\infty$. If*

$$\liminf_{r \to +\infty} \frac{Q(r) \text{vol} \, B_r}{r^2},$$

*then $M$ admits no conformal metric with scalar curvature $K(x)$.*

Now we will apply Theorem A' to the following geometric situation. Let $f : (M, <,>) \to \mathbb{R}^n$ be an isometric immersion with parallel mean curvature $H$ and arbitrary codimension. We assume that $M$ is oriented, and let $\gamma_f(M, <,>) \to G_m(\mathbb{R}^n)$ denote the Gauss map of $f$, where $G_m(\mathbb{R}^n)$ is the Grassmannian of oriented $m$-planes in $\mathbb{R}^n$. Identifying $m$-planes in $G_m(\mathbb{R}^n)$ with unit multivectors, we can define the cosine of the angle $\theta$ between $\gamma_f$ and the $m$-plane $V = V_1 \wedge V_2 \wedge ... \wedge V_m$ as

$$\cos \theta = <V_1 \wedge V_2 \wedge ... \wedge V_m, e_1 \wedge e_2 \wedge ... \wedge e_m>,$$

where $\{e_a\}$ is a local orthonormal frame along $f$, with $e_i$ tangent to $M$ for $i = 1, ..., m$, and $e_1 \wedge e_2 \wedge ... \wedge e_m$ giving the orientation, and $e_\alpha$ orthogonal to $M$ for $\alpha = m+1, ..., n$ (for short a local Darboux frame along $f$). Denoting by $S^{\binom{n}{m}-1}$ the unit sphere in $\Lambda^m(\mathbb{R}^n)$, we say that $\gamma_f(p)$ is contained in the open (respectively closed) spherical cap $S^{\binom{n}{m}-1}$ of radius $\theta_o$ $(0 \leq \theta_o \leq \pi)$ centered at $V_1 \wedge V_2 \wedge ... \wedge V_m$ if and only if $\cos(\theta) > \cos(\theta_o)$ (respectively $\cos(\theta) \geq \cos(\theta_o)$ ). Then we have the following theorem, which extends a result obtained in [RRS2, Theorem 1.2] under Ricci curvature assumptions.

**Theorem 1.4.** *Let $f : (M, <,>) \to \mathbb{R}^n$ be a complete, oriented, isometric immersion with parallel mean curvature $H$. Assume that the scalar curvature $s(x)$ of $M$ satisfies*

$$(1.2) \qquad s(x) \leq m^2|H|^2 - \frac{1}{Q(r(x))},$$

*where $Q(t)$ is a positive, non-decreasing function on $[0, +\infty)$ satisfying $Q(t) = o(t^2)$ as $t \to +\infty$. Suppose furthermore that*

$$(1.3) \qquad \liminf_{r \to +\infty} \frac{Q(r) \log \operatorname{vol} B_r}{r^2} < +\infty.$$

*Then $\gamma_f(M)$ is not contained in any closed spherical cap in $S^{\binom{n}{m}-1}$ of radius*

$$(1.4) \qquad \theta < \arccos\sqrt{\frac{2(n-m-1)}{3n-3m-2}}.$$

We conclude with an application to ordinary differential inequalities of the form

$$(1.5) \qquad \begin{cases} \left[\varphi(|\alpha'|)\operatorname{sgn}(\alpha')\right]' + \tau(r)\varphi(|\alpha'|)\operatorname{sgn}(\alpha') \geq b(r)f(\alpha) \\ \alpha(0) = \alpha_0, \quad \alpha'(0) = 0, \end{cases}$$

where $0 < b(r) \in C^0([0, +\infty))$, $\tau \in C^0((0, +\infty))$, $\tau(r) = O(1/r)$ as $r \to 0+$, $f \in C(\mathbb{R})$, and $\varphi$ satisfies the conditions listed in the Introduction.

These inequalities arise when one applies the $\varphi$-Laplacian on models to radial functions. Namely, given a positive function $g \in C^\infty([0, +\infty))$ such that

$$(1.6) \qquad g(0) = 0, \quad g'(0) = 1, \quad g^{(2k)}(0) = 0, \ k = 1, 2, \ldots,$$

let $<,>$ be the metric on $(0, +\infty) \times S^{m-1}$ defined by $<,> = dr^2 + g(r)^2 d\theta^2$. Then $<,>$ extends to a smooth metric on $\mathbb{R}^m$ for which

$$\operatorname{vol} B_r = c_m \int_0^r g^{m-1}(t)\,dt.$$

Note that if $\alpha \in C^2([0, +\infty)$ is such that $\alpha'(0) = 0$, and $u$ is the radial function on $(\mathbb{R}^m, <,>)$ defined by $u(x) = \alpha(r(x))$, then

$$\mathrm{d}iv\big(|\nabla u|^{-1}\varphi(|\nabla u|)\nabla u\big) = \left[\varphi(|\alpha'|)\operatorname{sgn}(\alpha')\right]' + (m-1)\frac{g'}{g}\varphi(|\alpha'|)\operatorname{sgn}(\alpha').$$

We refer to the paper by B. Franchi, E. Lanconelli and Serrin, [FLS], where existence and uniqueness of positive solutions is established under minimal conditions on $f$, in the case where the left hand side of (1.5) is the radial part of the $\varphi$-Laplacian on $\mathbb{R}^n$.

An application of Theorem A, yields the following

**Theorem 1.5.** *Let $\tau \in C^0((0,+\infty))$ satisfy*

$$(1.7) \qquad \tau(r) = \begin{cases} O(r^{-1}) & \text{as } r \to 0+ \\ O(r^\lambda) & \text{as } r \to +\infty, \end{cases}$$

*for some $-1 \leq \lambda \leq \delta$, and let $f(t) \geq 0$ for $t \geq 0$. Assume that $b(r) > 0$ on $[0,+\infty)$ and*

$$(1.8) \qquad b(r) \geq C_0 \begin{cases} \dfrac{\log r}{r^{1+\delta}} & \text{if } \lambda = -1 \\ \dfrac{1}{r^{\delta-\lambda}} & \text{if } -1 < \lambda \leq \delta, \end{cases}$$

*for some constant $C_0 > 0$ and sufficiently large $r$. Then (1.5) has no non-negative, bounded solution except, possibly, the solution $\alpha \equiv 0$. If we further assume that $f$ satisfies*

$$(1.9) \qquad \liminf \frac{f(t)}{t^\sigma} > 0$$

*for some $\sigma > 1$, then (1.5) has no non-negative solution, except, possibly, the solution $\alpha \equiv 0$.*

In the case of the Laplacian, a more refined version of Theorem 1.5 was obtained in [RRS1, Theorem 3.1] using entirely different methods. We note that there, the asymptotic lower bound for $b$ is better by log terms, and (1.9) is replaced by an integral condition on $f$. Further, $b$ is only required to be non-negative. However in order to deal with the case where $b$ vanishes, one uses methods based on conformal changes of metric that are unavailable in the more general setting covered by Theorem 1.5.

## References

[FLS]  B. Franchi, E. Lanconelli and j. Serrin, *Existence and uniqueness of non-negative solutions of quasilinear equations in $\mathbf{R}^n$*, Adv. Math. **118** (1996), no. 2, 177–243.

[G1]  A. Grigor'yan, *Bounded solutions of the Schrödinger equations on noncompact Riemannian manifolds* (Russian) Trudy Sem. Petrovsk. **14** (1989), 66–77, 265–266; translation in J. Soviet Math. **51** (1990), no. 3, 2340–2349.

[G2]   A. Grigor'yan, *Analytic and geometric background of recurrence and non-explosion of the Brownian motion on Riemannian manifolds*, Bull. Amer. Math. Soc. **36** (1999), 135–249.
[GK]   A. Grigor'yan and M. Kelbert, *On Hardy-Littlewood inequality for Brownian motion on Riemannian manifolds*, J. London Math. Soc. (2) **62** (2000), no. 2, 625–639.
[PeS]  L.A. Peletier and J. Serrin, *Gradient bounds and Liuoville theorems for quasilinear elliptic equations*, Ann. Scuola Norm. Sup. Pisa Cl. Sci. (4) **5** (1978) 65–104.
[PRS1] S. Pigola, M. Rigoli and A.G. Setti, *Maximum Principles and Singular Elliptic Inequalities*, J. Func. Anal., to appear.
[PRS2] S. Pigola, M. Rigoli and A.G. Setti, *A Remark on the Maximum Principle and Stochastic Completeness*, Proc. AMS, to appear.
[PRS3] S. Pigola, M. Rigoli and A.G. Setti, *Volume growth, "a priori" estimates and geometric applications*, preprint.
[PuS]  P. Pucci and J. Serrin, *A note on the strong maximum principle for elliptic differential inequalities*, J. Math. Pures Appl. **79** (2000), 57-71.
[PuSZ] P. Pucci, J. Serrin and H. Zou, *A strong maximum principle and a compact support principle for singular elliptic inequalities*, J. Math. Pures Appl. **78** (1999), 769–789.
[RRS1] A. Ratto, M. Rigoli and A.G. Setti, *On the Omori–Yau maximum principle and its applications to differential equations and geometry*, J. Func. Anal. **134** (1995), 486–510.
[RRS2] A. Ratto, M. Rigoli and A.G. Setti, *A uniqueness result in PDE's and parallel mean curvature immersions in Euclidean space*, Complex Variables **30** (1996), 221–233.
[RRV]  A. Ratto, M. Rigoli and L. Veron, *Scalar curvature and conformal deformations of noncompact manifolds*, Math. Z. **225** (1997), 395–426.
[RS]   M. Rigoli and A.G. Setti, *Liouville-type theorems for $\varphi$-subharmonic functions*, Rev. Mat. Iberoamericana, to appear.
[S1]   J. Serrin, *Entire solutions of nonlinear Poisson equations*, Proc. London Math. Soc. **24** (1972), 384-366.
[S2]   J. Serrin, *On the strong maximum principle for quasilinear second order differential inequalities*, J. Func. Anal. **5** (1970), 184–193.

Dipartimento di Matematica, Università di Milano, via Saldini 50, I-20133 Milano, ITALY
  *E-mail address*: `pigola@mat.unimi.it`

Dipartimento di Scienze Chimiche, Fisiche e Matematiche, Facoltà di Scienze, Università dell'Insubria - Como, via Valleggio 11, I-22100 Como, ITALY
  *E-mail address*: `rigoli@matapp.unimib.it`

Dipartimento di Scienze Chimiche, Fisiche e Matematiche, Facoltà di Scienze, Università dell'Insubria - Como, via Valleggio 11, I-22100 Como, ITALY
  *E-mail address*: `setti@uninsubria.it`